\documentclass[a4paper,10pt]{amsart}      
\usepackage{amssymb,amsthm, amsmath, amsfonts} 
\usepackage{graphics}                 
\usepackage{color}                    
\usepackage{hyperref}                 
\usepackage[all]{xy}

\newtheorem{theorem}{Theorem}[section]
\newtheorem{lemma}[theorem]{Lemma}
\newtheorem{proposition}[theorem]{Proposition}

\newtheorem{corollary}[theorem]{Corollary}
\newtheorem{definition}[theorem]{Definition}
\newtheorem{example}[theorem]{Example}
\newtheorem{remark}[theorem]{Remark}
\numberwithin{equation}{section}
\def\rad{\mathop{\hbox{rad}}}
\def\Ob{\mathop{\hbox{Ob}}}

\def\dim{\mathop{\hbox{dim}}}
\def\top{\mathop{\hbox{Top}}}

\title{Extension Algebras of Standard Modules}
\author{Liping Li}
\address{School of Mathematics, University of Minnesota, MN, 55455, USA}
\email{lixxx480@math.umn.edu}
\thanks{The author would like to thank his thesis advisor, Professor Peter Webb, for proposing to study the the extension algebras of standard modules, and carefully checking the manuscript. He also thanks Prof. Mazorchuk for pointing out many related works, which are unknown to the author before.}

\begin{document}
\begin{abstract}
Let $A$ be a basic finite-dimensional $k$-algebra standardly stratified for a partial order $\leqslant$ and $\Delta$ be the direct sum of all standard modules. In this paper we study the extension algebra $\Gamma = \text{Ext} _A^{\ast} (\Delta, \Delta)$ of standard modules, characterize the stratification property of $\Gamma$ for $\leqslant$ and $\leqslant ^{op}$, and obtain a sufficient condition for $\Gamma$ to be a generalized Koszul algebra (in a sense which we define).
\end{abstract}
\maketitle

Let $A$ be a basic finite-dimensional $k$-algebra standardly stratified with respect to a poset $(\Lambda, \leqslant)$ indexing all simple modules (up to isomorphism), $\Delta$ be the direct sum of all standard modules, and $\mathcal{F} (\Delta)$ be the category of finitely generated $A$-modules with $\Delta$-filtrations. That is, for each $M \in \mathcal{F} (\Delta)$, there is a chain $0 = M_0 \subseteq M_1 \subseteq \ldots \subseteq M_n = M$ such that $M_i / M_{i-1}$ is isomorphic to an indecomposable summand of $\Delta$, $1 \leqslant i \leqslant n$. Since standard modules of $A$ are relative simple in $\mathcal{F} (\Delta)$, we are motivated to exploit the extension algebra $\Gamma = \text{Ext} _A^{\ast} (\Delta, \Delta)$ of standard modules. These extension algebras were studied in \cite{Abe, Drozd, Klamt, Mazorchuk1, Miemietz}. In this paper, we are interested in the stratification property of $\Gamma$ with respect to $(\Lambda, \leqslant)$ and $(\Lambda, \leqslant ^{op})$, and its Koszul property since $\Gamma$ has a natural grading. A particular question is that in which case it is a \textit{generalized Koszul algebra}, i.e., $\Gamma_0$ has a linear projective resolution.

By Gabriel's construction, we associate a locally finite $k$-linear category $\mathcal{E}$ to the extension algebra $\Gamma$ such that the category $\Gamma$-mod of finitely generated left $\Gamma$-modules is equivalent to the category of finitely generated $k$-linear representations of $\mathcal{E}$. We show that the category $\mathcal{E}$ is a \textit{directed category} with respect to $\leqslant$. That is, the morphism space $\mathcal{E} (x,y) = 0$ whenever $x \nleqslant y$. With this terminology, we have:

\begin{theorem}
If $A$ is standardly stratified for $(\Lambda, \leqslant)$, then $\mathcal{E}$ is a directed category with respect to $\leqslant$ and is standardly stratified for $\leqslant ^{op}$. Moreover, $\mathcal{E}$ is standardly stratified for $\leqslant$ if and only if for all $\lambda, \mu \in \Lambda$ and $s \geqslant 0$, Ext$ _A^s (\Delta_{\lambda}, \Delta_{\mu})$ is a projective End$ _A (\Delta_{\mu})$-module.
\end{theorem}

In particular, if $A$ is a quasi-hereditary algebra, then $\Gamma$ is quasi-hereditary with respect to both $\leqslant$ and $\leqslant ^{op}$. We also generalize the above theorem to abstract \textit{stratifying systems} and \textit{Ext-projective stratifying systems} (EPSS) described in \cite{Erdmann, Marcos1, Marcos2, Webb}.

In the case that the standardly stratified algebra $A$ is a graded algebra and $A_0$ is semisimple, its Koszul duality has been studied in \cite{Agoston1, Agoston2, Mazorchuk2, Mazorchuk3, Mazorchuk4}. Since the extension algebra $\Gamma = \text{Ext} _A^{\ast} (\Delta, \Delta)$ has a natural grading, and $\Gamma_0 = \text{End} _A(\Delta)$ in general is not semisimple, we study the general Koszul property of $\Gamma$ by using some generalized Koszul theories developed in \cite{Green, Li1, Madsen2, Madsen3}, where the degree 0 parts of graded algebras are not required to be semisimple.

Take a fixed EPSS $(\underline {\Theta}, \underline{Q})$. As an analogue to linear modules of graded algebras, we define \textit{linearly filtered modules} in this system. With this terminology, a sufficient condition can be obtained for $\Gamma$ to be a generalized Koszul algebra.

\begin{theorem}
Let $(\underline {\Theta}, \underline{Q})$ be an EPSS indexed by a finite poset $(\Lambda, \leqslant)$ such that Ext$_A ^i (Q, \Theta) =0$ for all $i \geqslant 1$ and Hom$ _A (Q, \Theta) \cong \text{Hom} _A (\Theta, \Theta)$. Suppose that all $\Theta_{\lambda}$ are linearly filtered for $\lambda \in \Lambda$. If $M \in \mathcal{F} (\Theta)$ is linearly filtered, then the graded $\Gamma$-module $\text{Ext} _A^{\ast} (M, \Theta)$ has a linear projective resolution. In particular, $\Gamma = \text{Ext} _A^{\ast} (\Theta, \Theta)$ is a generalized Koszul algebra.
\end{theorem}

The paper is organized as follows: in section 1 we characterize the stratification property of $\Gamma$; in section 2 we define linearly filtered modules, study their basic properties and prove the second theorem.

Throughout this paper $A$ is a finite-dimensional basic associative $k$-algebra with identity 1, where $k$ is algebraically closed. We only consider finitely generated modules and denote by $A$-mod the category of finitely generated left $A$-modules. Maps and morphisms are composed from right to left.

\section{Stratification property of extension algebras}

Let $(\Lambda, \leqslant)$ be a finite preordered set parameterizing all simple $A$-modules $S_{\lambda}$ (up to isomorphism). This preordered set also parameterizes all indecomposable projective $A$-modules $P_{\lambda}$ (up to isomorphism). According to \cite{Cline}, the algebra $A$ is standardly-stratified with respect to $(\Lambda, \leqslant)$ if there exist modules $\Delta_{\lambda}$, $\lambda \in \Lambda$, such that the following conditions hold:
\begin{enumerate}
\item the composition factor multiplicity $[\Delta_{\lambda} : S_{\mu}] = 0$ whenever $\mu \nleqslant \lambda$;
and
\item for every $\lambda \in \Lambda$ there is a short exact sequence $0 \rightarrow K_{\lambda} \rightarrow P_{\lambda} \rightarrow \Delta_{\lambda} \rightarrow 0$ such that $K_{\lambda}$ has a filtration with factors $\Delta_{\mu}$ where $\mu > \lambda$.
\end{enumerate}
In some literatures the preordered set $(\Lambda, \leqslant)$ is supposed to be a poset (\cite{Dlab}) or even a linearly ordered set (\cite{Agoston1, Agoston2}). Algebras standardly stratified in this sense are called \textit{strongly standardly stratified} (\cite{Frisk1, Frisk2}). In this paper $(\Lambda, \leqslant)$ is supposed to be a poset.

If $A$ is standardly stratified, then \textit{standard modules} can be defined as \footnote{In \cite{Agoston1, Agoston2} standard modules are defined as $\Delta_{\lambda} = P_{\lambda} / \sum _{\mu > \lambda} \text{tr} _{P_{\mu}} (P_{\lambda})$. Note that in their setup $\leqslant$ is a linear order, so this description of standard modules coincides with ours.}:
\begin{equation*}
\Delta_{\lambda} = P_{\lambda} / \sum _{\mu \nleqslant \lambda} \text{tr} _{P_{\mu}} (P_{\lambda}),
\end{equation*}
where tr$_{P_{\mu}} (P_{\lambda})$ is the trace of $P_{\mu}$ in $P _{\lambda}$. See \cite{Dlab, Webb} for more details. Let $\Delta$ be the direct sum of all standard modules and $\mathcal{F} (\Delta)$ be the full subcategory of $A$-mod such that each object in $\mathcal{F} (\Delta)$ has a filtration by standard modules. Clearly, since $A$ is standardly stratified for $\leqslant$, $_AA \in \mathcal{F} (\Delta)$, or equivalently, every indecomposable projective $A$-module has a filtration by standard modules.

Throughout this section we suppose that $A$ is standardly stratified with respect to $\leqslant$ if it is not specified. We also remind the reader that $\mathcal{F} (\Delta)$ is closed under extensions, kernels of epimorphisms, and direct summands, but is not closed under cokernels of monomorphisms.

Given $M \in \mathcal{F} (\Delta)$ and a fixed filtration $0 = M_0 \subseteq M_1 \subseteq \ldots \subseteq M_n =M$, we define the \textit{filtration multiplicity} $m_{\lambda} = [ M: \Delta_{\lambda}]$ to be the number of factors isomorphic to $\Delta_{\lambda}$ in this filtration. By Lemma 1.4 of \cite{Erdmann}, The filtration multiplicities defined above are independent of the choice of a particular filtration. Moreover, since each standard module has finite projective dimension, we deduce that every $A$-module contained in $\mathcal{F} (\Delta)$ has finite projective dimension. Therefore, the extension algebra $\Gamma = \text{Ext} _A^{\ast} (\Delta, \Delta)$ is finite-dimensional.

\begin{lemma}
Let $\Delta_{\lambda}$, $\Delta_{\mu}$ be standard modules. Then Ext$ _A^n (\Delta _{\lambda}, \Delta_{\mu}) =0$ if $\lambda \nleqslant \mu$ for all $n \geqslant 0$.
\end{lemma}

\begin{proof}
First, we claim $[\Omega^i (\Delta_{\lambda}) : \Delta_{\nu}] = 0$ whenever $\lambda \nleqslant \nu$ for all $i \geqslant 0$, where $\Omega$ is the Heller operator. Indeed, for $i =0$ the conclusion holds clearly. Suppose that it is true for all $i \leqslant n$ and consider $\Omega ^{n+1} (\Delta _{\lambda})$. We have the following exact sequence:
\begin{equation*}
\xymatrix {0 \ar[r] & \Omega ^{n+1} (\Delta _{\lambda}) \ar[r] & P \ar[r] & \Omega^n (\Delta _{\lambda}) \ar[r] & 0.}
\end{equation*}
By the induction hypothesis, $[\Omega ^n (\Delta _{\lambda}) : \Delta_{\nu}] = 0$ whenever $\lambda \nleqslant \nu$. Therefore, $[P : \Delta_{\nu}] = 0$ whenever $\lambda \nleqslant \nu$, and hence $[\Omega ^{n+1} (\Delta _{\lambda}) : \Delta_{\nu}] = 0$ whenever $\lambda \nleqslant \nu$. The claim is proved by induction.

The above short exact sequence induces a surjection Hom$ _A (\Omega ^n (\Delta _{\lambda}), \Delta_{\mu}) \rightarrow \text{Ext} ^n_A (\Delta _{\lambda}, \Delta_{\mu})$. Thus it suffices to show Hom$ _A (\Omega ^n (\Delta _{\lambda}), \Delta_{\mu}) =0$ for all $n \geqslant 0$ if $\lambda \nleqslant \mu$. By the above claim, all filtration factors $\Delta_{\nu}$ of $\Omega ^n (\Delta _{\lambda})$ satisfy $\nu \geqslant \lambda$, and hence $\nu \nleqslant \mu$. But Hom$ _{A} (\Delta_{\nu}, \Delta_{\mu}) = 0$ whenever $\nu \nleqslant \mu$. The conclusion follows.
\end{proof}

Gabriel's construction gives rise to a bijective correspondence between finite-dimensional algebras and locally finite $k$-linear categories with finitely many objects. Explicitly, To each finite-dimensional $k$-algebra $A$ with a chosen set of orthogonal primitive idempotents $\{ e_{\lambda} \} _{\lambda \in \Lambda}$ satisfying $\sum _{\lambda \in \Lambda} e_{\lambda} =1$ we define a $k$-linear category $\mathcal{A}$ with $\Ob \mathcal{A} = \{ e_{\lambda} \} _{\lambda \in \Lambda}$ and $\mathcal{A} (e_{\lambda}, e_{\mu}) = e_{\mu} A e_{\lambda} \cong \text{Hom} _{A} (A e_{\mu}, A e_{\lambda})$. Conversely, given a locally finite $k$-linear category $\mathcal{A}$ with finitely many objects, we define $A = \bigoplus _{x, y \in \text{Ob } \mathcal{A}} \mathcal{A} (x,y)$, and the multiplication in $A$ is induces by the composition of morphisms in $\mathcal{A}$ in an obvious way. Clearly, $A$-mod is Morita equivalent to the category of all finite-dimensional $k$-linear representations of $\mathcal{A}$. We then call $\mathcal{A}$ the \textit{associated category} of $A$ and $A$ the \textit{associated algebra} of $\mathcal{A}$. The category $\mathcal{A}$ is called \textit{directed} if there is a partial order $\preccurlyeq$ on $\Ob \mathcal{A}$ such that $\mathcal{A} (x, y) = 0$ unless $x \preccurlyeq y$.

Now let $\Gamma = \text{Ext} _A^{\ast} (\Delta, \Delta)$. This is a graded finite-dimensional algebra equipped with a natural grading. In particular, $\Gamma_0 = \text{End} _A (\Delta)$. For each $\lambda \in \Lambda$, $\Delta_{\lambda}$ is an indecomposable $A$-module. Therefore, up to isomorphism, the indecomposable projective $\Gamma$-modules are exactly those $\text{Ext} _A^{\ast} (\Delta_{\lambda}, \Delta)$, $\lambda \in \Lambda$.

The associated $k$-linear category $\mathcal{E}$ of $\Gamma$ has the following structure: $\Ob \mathcal{E} = \{ \Delta _{\lambda} \} _{\lambda \in \Lambda}$; the morphism space $\mathcal{E} (\Delta _{\lambda}, \Delta_{\mu}) = \text{Ext} _A^{\ast} (\Delta_{\lambda}, \Delta_{\mu})$. The partial order $\leqslant$ induces a partial order on $\Ob \mathcal{E}$ which we still denote by $\leqslant$, namely, $\Delta _{\lambda} \leqslant \Delta_{\mu}$ if and only if $\lambda \leqslant \mu$.

\begin{proposition}
The associated category $\mathcal{E}$ of $\Gamma$ is directed with respect to $\leqslant$. In particular, $\Gamma$ is standardly stratified with respect to $\leqslant ^{op}$ and all standard modules are projective.
\end{proposition}

\begin{proof}
The first statement follows from the previous lemma. The second statement is also clear. Indeed, since $\Gamma$ is directed with respect to $\leqslant$, $e_{\mu} \Gamma e_{\lambda} \cong \text{Hom} _{\Gamma} (Q_{\mu}, Q_{\lambda}) = 0$ if $\mu \ngeqslant \lambda$, where $Q_{\mu}, Q_{\lambda}$ are projective $\Gamma$-modules. Thus tr$ _{Q_{\mu}} (Q_{\lambda}) = 0$ whenever $\mu \ngeqslant \lambda$, or equivalently, tr$ _{Q_{\mu}} (Q_{\lambda}) = 0$ whenever $\mu \nleqslant ^{op} \lambda$. Therefore, all standard modules with respect to $\leqslant ^{op}$ are projective.
\end{proof}

The following proposition characterizes the stratification property of a $k$-linear category directed with respect to $\leqslant$.

\begin{proposition}
Let $\mathcal{C}$ be a locally finite $k$-linear category directed with respect to a partial order $\leqslant$ on $\Ob \mathcal{C}$. Then it is stratified for this order. The standard modules are isomorphic to indecomposable summands of $\bigoplus _{x \in \text{Ob } \mathcal{C}} \mathcal{C} (x, x)$. Moreover, this stratification is standard if and only if for each pair of objects $x, y \in \Ob \mathcal{C}$, $\mathcal{C} (x, y)$ is a projective $\mathcal{C} (y,y)$-module.
\end{proposition}

\begin{proof}
This is just a collection of results in \cite{Li1}. The first statement is Corollary 5.4; the second statement comes from Proposition 5.5; and the last statement is Theorem 5.7.
\end{proof}

Now we restate and prove the first theorem.

\begin{theorem}
If $A$ is standardly stratified for $(\Lambda, \leqslant)$, then $\mathcal{E}$ is a directed category with respect to $\leqslant$ and is standardly stratified for $\leqslant ^{op}$. Moreover, $\mathcal{E}$ is standardly stratified for $\leqslant$ if and only if for all $\lambda, \mu \in \Lambda$ and $s \geqslant 0$, Ext$ _A^s (\Delta_{\lambda}, \Delta_{\mu})$ is a projective End$ _A (\Delta_{\mu})$-module.
\end{theorem}

\begin{proof}
The first statement follows from Proposition 1.2 and the second statement follows from Proposition 1.3.
\end{proof}

In the case that $A$ is quasi-hereditary, we have:

\begin{corollary}
If $A$ is a quasi-hereditary algebra with respect to $\leqslant$, then $\Gamma$ is quasi-hereditary with respect to both $\leqslant$ and $\leqslant ^{op}$.
\end{corollary}

\begin{proof}
We have showed that $\Gamma$ is standardly stratified with respect to $\leqslant ^{op}$ and the corresponding standard modules $_{\Gamma} \Delta _{\lambda} \cong \Gamma 1_{\lambda}$ for $\lambda \in \Lambda$. Therefore,
\begin{equation*}
\text{End} _{\Gamma} (_{\Gamma} \Delta _{\lambda}) = \text{End} _{\Gamma} (\Gamma 1_{\lambda}) \cong 1_{\lambda} \Gamma 1_{\lambda} = \text{Ext} _A^{\ast} (\Delta_{\lambda}, \Delta_{\lambda}) = \text{End}_A (\Delta_{\lambda}) \cong k
\end{equation*}
since $A$ is quasi-hereditary. So $\Gamma$ is also quasi-hereditary with respect to $\leqslant^{op}$.

Now consider the stratification property of $\Gamma$ with respect to $\leqslant$. The associated category $\mathcal{E}$ is directed with respect to $\leqslant$. Since $\text{Ext} _A^{\ast} (\Delta_{\mu}, \Delta_{\mu}) = \text{End} _A (\Delta_{\mu}) \cong k$ for all $\mu \in \Lambda$, $\mathcal{E} (\Delta_{\lambda}, \Delta_{\mu}) = \text{Ext} _A^{\ast} (\Delta_{\lambda}, \Delta_{\mu})$ is a projective $k$-module for each pair $\lambda, \mu \in \Lambda$. Therefore, $\mathcal{E}$ is standardly stratified for $\leqslant$ by the previous theorem. Moreover, by Proposition 1.3, the standard modules of $\mathcal{E}$ (or the standard modules of $\Gamma$) are precisely indecomposable summands of $\bigoplus _{\lambda \in \Lambda} \text{Ext} _A^{\ast} (\Delta_{\lambda}, \Delta_{\lambda}) \cong \bigoplus _{\lambda \in \Lambda} k_{\lambda}$. Clearly, for $\lambda \in \Lambda$, End$ _{\Gamma} (k_{\lambda}, k_{\lambda}) \cong k$, so $\Gamma$ is quasi-hereditary with respect to $\leqslant$.
\end{proof}

The following example from 8.2 in \cite{Frisk2} illustrates why we should assume that $\leqslant$ is a partial order rather than a preorder. Indeed, in a preordered set $(\Lambda, \leqslant)$ we cannot deduce $x = y$ if $x \leqslant y$ and $y \leqslant x$.

\begin{example}
Let $A$ be the path algebra of the following quiver with relations $\alpha_1 \beta_1 = \alpha_2 \beta_2 = \alpha_2 \alpha_1 = \beta_1 \beta_2 =0$. Define a preorder $\leqslant$ by letting $x \leqslant y < z$ and $y \leqslant x < z$.
\begin{equation*}
\xymatrix{x \ar@/^/[r] ^{\alpha_1} & y \ar@/^/[r] ^{\alpha_2} \ar@/^/[l] ^{\beta_1} & z \ar@/^/[l] ^{\beta_2}}.
\end{equation*}
Projective modules and standard modules are described as follows:
\begin{equation*}
P_x \cong \Delta_x = \begin{matrix} x \\ y \\ x \end{matrix} \qquad P_y = \begin{matrix} & y & \\ x & & z \\ & & y \end{matrix} \qquad \Delta_y = \begin{matrix} y \\ x \end{matrix} \qquad P_z \cong \Delta_z = \begin{matrix} z \\ y \end{matrix}
\end{equation*}
Then the associated category $\mathcal{E}$ of $\Gamma = \text{Ext} _A^{\ast} (\Delta, \Delta)$ is not a directed category since both Hom$_A (\Delta_x, \Delta_y)$ and Hom$ _A (\Delta_y, \Delta_x)$ are nonzero.
\end{example}

Now we generalize the above results to \textit{Ext-Projective Stratifying Systems} (EPSS). From now on the algebra $A$ is finite-dimensional and basic, but we do not assume that it is standardly stratified for some partial order, as we did before. The EPSS we describe in this paper is indexed by a finite poset $(\Lambda, \leqslant)$ rather than a linearly ordered set as in \cite{Marcos1, Marcos2}. However, this difference is not essential and all properties described in \cite{Marcos1, Marcos2} can be applied to our situation with suitable modifications.

\begin{definition}
(Definition 2.1 in \cite{Marcos2}) Let $\underline{\Theta} = \{ \Theta_{\lambda} \}_{\lambda \in \Lambda}$ be a set of nonzero $A$-modules and $\underline{Q} = \{ Q_{\lambda} \} _{\lambda \in \Lambda}$ be a set of indecomposable $A$-modules, both of which are indexed by a finite poset $(\Lambda, \leqslant)$. We call $(\underline {\Theta}, \underline{Q})$ an EPSS if the following conditions are satisfied:
\begin{enumerate}
\item Hom$ _A (\Theta_{\lambda}, \Theta_{\mu}) =0$ if $\lambda \nleqslant \mu$;
\item for each $\lambda \in \Lambda$, there is an exact sequence $0 \rightarrow K_{\lambda} \rightarrow Q_{\lambda} \rightarrow \Theta_{\lambda} \rightarrow 0$ such that $K_{\lambda}$ has a filtration only with factors isomorphic to $\Theta_{\mu}$ satisfying $\mu > \lambda$;
\item for every $A$-module $M \in \mathcal{F} (\underline {\Theta})$ and $\lambda \in \Lambda$, Ext$ _A^1 (Q_{\lambda}, M) =0$.
\end{enumerate}
\end{definition}

We denote $\Theta$ and $Q$ the direct sums of all $\Theta _{\lambda}$'s and $Q_{\lambda}$'s respectively, $\lambda \in \Lambda$.

Given an EPSS $(\underline {\Theta}, \underline{Q})$ indexed by $(\Lambda, \leqslant)$, $(\underline {\Theta}, \leqslant)$ is a \textit{stratifying system} (SS): Hom$ _A (\Theta_{\lambda}, \Theta_{\mu}) =0$ if $\lambda \nleqslant \mu$, and Ext$ _A^1 (\Theta_{\lambda}, \Theta_{\mu}) =0$ if $\lambda \nless \mu$. Conversely, given a stratifying system $(\underline {\Theta}, \leqslant)$, we can construct an EPSS $(\underline {\Theta}, \underline{Q})$ unique up to isomorphism. See \cite{Marcos2} for more details. Moreover, as described in \cite{Marcos2}, the algebra $B = \text{End} _A (Q) ^{op}$ is standardly stratified, and the functor $e_Q = \text{Hom} _A (Q, -)$ gives an equivalence of exact categories between $\mathcal{F} (\Theta)$ and $\mathcal{F} (_B \Delta)$.

To study the extension algebra $\Gamma = \text{Ext} _A^{\ast} (\Theta, \Theta)$, one may want to use projective resolutions of $\Theta$. However, different from the situation of standardly stratified algebras, the regular module $_AA$ in general might not be contained in $\mathcal{F} (\Theta)$. If we suppose that $_AA$ is contained in $\mathcal{F} (\Theta)$ (in this case the stratifying system $(\underline {\Theta}, \leqslant)$ is said to be \textit{standard}) and $\mathcal{F} (\Theta)$ is closed under the kernels of surjections, then by Theorem 2.6 in \cite{Marcos1} $A$ is standardly stratified for $\leqslant$ and those $\Theta_{\lambda}$'s coincide with standard modules of $A$. This situation has been completely discussed previously. Alternately, we use the \textit{relative projective resolutions} whose existence is guaranteed by the following proposition.

\begin{proposition}
(Corollary 2.11 in \cite{Marcos2}) Let $(\underline {\Theta}, \underline{Q})$ be an EPSS indexed by a finite poset $(\Lambda, \leqslant)$. Then for each $M \in \mathcal{F} (\Theta)$, there is a finite resolution
\begin{equation*}
\xymatrix {0 \ar[r] & Q^d \ar[r] & \ldots \ar[r] & Q^0 \ar[r] & M \ar[r] & 0}
\end{equation*}
such that each kernel is contained in $\mathcal{F} (\Theta)$, where $0 \neq Q^i \in \text{add} (Q)$ for $0 \leqslant i \leqslant d$.
\end{proposition}

The number $d$ in this resolution is called the \textit{relative projective dimension} of $M$.

\begin{proposition}
Let $(\underline {\Theta}, \underline{Q})$ be an EPSS indexed by a finite poset $(\Lambda, \leqslant)$ and $d$ be the relative projective dimension of $\Theta$. If Ext$ _A^s (Q, \Theta) = 0$ for all $s \geqslant 1$, then for $M, N \in \mathcal{F} (\Theta)$ and $s > d$, Ext$ _A^s (M, N) =0$.
\end{proposition}

\begin{proof}
Since both $M$ and $N$ are contained in $\mathcal{F} (\Theta)$, it is enough to show that Ext$ _A^s (\Theta, \Theta) = 0$ for all $s > d$. If $d =0$, then $Q = \Theta$ and the conclusion holds trivially. So we suppose $d \geqslant 1$. Applying the functor Hom$ _A (-, \Theta)$ to the exact sequence
\begin{equation*}
\xymatrix {0 \ar[r] & K_1 \ar[r] & Q \ar[r] & \Theta \ar[r] & 0}
\end{equation*}
we get a long exact sequence. In particular, from the segment
\begin{equation*}
\xymatrix {\text{Ext} _A^{s-1} (Q, \Theta) \ar[r] & \text{Ext} _A^{s-1} (K_1, \Theta) \ar[r] & \text{Ext} _A^s (\Theta, \Theta) \ar[r] & \text{Ext} _A^s (Q, \Theta)}
\end{equation*}
of this long exact sequence we deduce that $\text{Ext} _A^s (\Theta, \Theta) \cong \text{Ext} _A^{s-1} (K_1, \Theta)$ since the first and last terms are 0. Now applying Hom$ _A (-, \Theta)$ to the exact sequence
\begin{equation*}
\xymatrix {0 \ar[r] & K_2 \ar[r] & Q^1 \ar[r] & K_1 \ar[r] & 0}
\end{equation*}
we get $\text{Ext} _A^{s-1} (K_1, \Theta) \cong \text{Ext} _A^{s-2} (K_2, \Theta)$. Thus $\text{Ext} _A^s (\Theta, \Theta) \cong \text{Ext} _A^{s-d} (K_d, \Theta)$ by induction. But $K_d \cong Q^d \in \text{add} (Q)$. The conclusion follows.
\end{proof}

Thus $\Gamma = \text{Ext} _A^{\ast} (\Theta, \Theta)$ is a finite-dimensional algebra under the given assumption.

There is a natural partition on the finite poset $(\Lambda, \leqslant)$ as follows: let $\Lambda_1$ be the subset of all minimal elements in $\Lambda$, $\Lambda_2$ be the subset of all minimal elements in $\Lambda \setminus \Lambda_1$, and so on. Then $\Lambda = \sqcup_{i \geqslant 1} \Lambda_i$. With this partition, we can introduce a \textit{height function} $h: \Lambda \rightarrow \mathbb{N}$ in the following way: for $\lambda \in \Lambda_i \subseteq \Lambda$, $i \geqslant 1$, we define $h(\lambda) = i$.

For each $M \in \mathcal{F} (\Theta)$, we define supp$(M)$ to be the set of elements $\lambda \in \Lambda$ such that $M$ has a $\Theta$-filtration in which there is a factor isomorphic to $\Theta_{\lambda}$. For example, supp$ (\Theta_{\lambda}) = \{ \lambda \}$. By Lemma 2.6 in \cite{Marcos2}, the multiplicities of factors of $M$ is independent of the choice of a particular $\Theta$-filtration. Therefore, supp$(M)$ is well defined. We also define $\min (M) = \min (\{ h(\lambda) \mid \lambda \in \text{supp} (M) \})$.

\begin{lemma}
Let $(\underline {\Theta}, \underline{Q})$ be an EPSS indexed by a finite poset $(\Lambda, \leqslant)$. For each $M \in \mathcal{F} (\Theta)$, there is an exact sequence $0 \rightarrow K_1 \rightarrow Q^0 \rightarrow M$ such that $K_1 \in \mathcal{F} (\Theta)$ and $\min (K_1) > \min (M)$, where $Q^0 \in \text{add} (Q)$.
\end{lemma}

\begin{proof}
This is Proposition 2.10 in \cite{Marcos2} which deals with the special case that $\Lambda$ is a linearly ordered set. The general case can be proved similarly by observing the fact that Ext$ _A^1 (\Theta_{\lambda}, \Theta_{\mu}) = 0$ if $h (\lambda) = h (\mu)$.
\end{proof}

By this lemma, the relative projective dimension of every $M \in \mathcal{F} (\Theta)$ cannot exceed the length of the longest chain in $\Lambda$.

As before, we let $\mathcal{E}$ be the $k$-linear category associated to $\Gamma = \text{Ext} _A^{\ast} (\Theta, \Theta)$.

\begin{theorem}
Let $(\underline {\Theta}, \underline{Q})$ be an EPSS indexed by a finite poset $(\Lambda, \leqslant)$. such that Ext$_A ^i (Q, \Theta) =0$ for all $i \geqslant 1$. Then $\mathcal{E}$ is a directed category with respect to $\leqslant$ and is standardly stratified for $\leqslant ^{op}$. Moreover, it is standardly stratified for $\leqslant$ if and only if for all $s \geqslant 0$, Ext$ _A^s (\Theta_{\lambda}, \Theta_{\mu})$ is a projective End$ _A (\Theta_{\mu})$-module, $\lambda, \mu \in \Lambda$.
\end{theorem}

\begin{proof}
We only need to show that $\mathcal{E}$ is a directed category with respect to $\leqslant$ since the other statements can be proved as in Theorem 1.4. We know Hom$ _A (\Theta_{\lambda}, \Theta_{\mu}) = 0$ if $\lambda \nleqslant \mu$ and Ext$ _A^1 (\Theta_{\lambda}, \Theta_{\mu}) = 0$ for all $\lambda \nless \mu$. Therefore, it suffices to show that for all $s \geqslant 2$, Ext$ _A^s (\Theta_{\lambda}, \Theta_{\mu}) = 0$ if $\lambda \nless \mu$.

By Proposition 1.8 and Lemma 1.10, $\Theta_{\lambda}$ has a relative projective resolution
\begin{equation*}
\xymatrix {0 \ar[r] & Q^d \ar[r] ^{f_d} & \ldots \ar[r] ^{f_1} & Q^0 \ar[r] ^{f_0} & \Theta_{\lambda} \ar[r] & 0}
\end{equation*}
such that for each map $f_t$, min($K_t) > \text{min} (K_{t-1})$, where $K_t = \text{Ker} (f_t)$ and $1 \leqslant t \leqslant d$. By Proposition 1.9, Ext$ _A^s (\Theta_{\lambda}, \Theta_{\mu}) = 0$ if $s > d$; if $2 \leqslant s \leqslant d$, we have Ext$ _A^s (\Theta_{\lambda}, \Theta_{\mu}) \cong \text{Ext} _A^1 (K_{s-1}, \Theta_{\mu})$. But we have chosen
\begin{equation*}
\min (K_{s-1}) > \min (K_{s-2}) > \ldots > \min (\Theta_{\lambda}) = h(\lambda) \geqslant h(\mu).
\end{equation*}
Thus each factor $\Theta_{\nu}$ appearing in a $\Theta$-filtration of $K_{s-1}$ satisfies $h(\nu) > h(\mu)$, and hence $\nu \nleqslant \mu$. Since Ext$ _A^1 (\Theta_{\nu}, \Theta_{\mu}) =0$ for all $\nu \nleqslant \mu$, we deduce
\begin{equation*}
\text{Ext} _A^s (\Theta_{\lambda}, \Theta_{\mu}) \cong \text{Ext} _A^1 (K_{s-1}, \Theta_{\mu}) =0.
\end{equation*}
This finishes the proof.
\end{proof}

The following corollary is a generalization of Corollary 1.5.

\begin{corollary}
Let $(\underline {\Theta}, \underline{Q})$ be an EPSS indexed by a finite poset $(\Lambda, \leqslant)$. If for all $s \geqslant 1$ and $\lambda \in \Lambda$ we have Ext$ _A ^s (Q, \Theta) = 0$ and End$_A (\Theta_{\lambda}, \Theta_{\lambda}) \cong k$, then $\Gamma$ is quasi-hereditary with respect to both $\leqslant$ and $\leqslant ^{op}$.
\end{corollary}

\begin{proof}
This can be proved as Corollary 1.5.
\end{proof}

\section{Koszul Property of Extension Algebras}

There is a well known duality related to the extension algebras: the Koszul duality. Explicitly, if $A$ is a graded Koszul algebra with $A_0$ being a semisimple algebra, then $B = \text{Ext} _A^{\ast} (A_0, A_0)$ is a Koszul algebra, too. Moreover, the functor $\text{Ext} _A^{\ast} (-, A_0)$ gives an equivalence between the category of linear $A$-modules and the category of linear $B$-modules.\footnote{In \cite{Li1} we generalized these results to the situation that $A_0$ is a self-injective algebra.} However, even if $A$ is quasi-hereditary with respect to a partial order $\leqslant$, $B$ might not be quasi-hereditary with respect to $\leqslant$ or $\leqslant ^{op}$. This problem has been considered in \cite{Agoston1, Mazorchuk2}.

On the other hand, if $A$ is quasi-hereditary with respect to $\leqslant$, we have showed that the extension algebra $\Gamma = \text{Ext} _A^{\ast} (\Delta, \Delta)$ is quasi-hereditary with respect to both $\leqslant$ and $\leqslant ^{op}$. But $\Gamma$ is in general not a Koszul algebra in a sense which we define later. In this section we want to get a sufficient condition for $\Gamma$ to be a generalized Koszul algebra.

We work in the context of EPSS described in last section. Let $(\underline {\Theta}, \underline{Q})$ be an EPSS indexed by a finite poset $(\Lambda, \leqslant)$; $Q = \bigoplus _{\lambda \in \Lambda} Q_{\lambda}$ and $\Theta = \bigoplus _{\lambda \in \Lambda} \Theta_{\lambda}$. \textbf{We insist the following conditions: Ext$ \mathbf{_A^s (Q, \Theta) =0}$ for all $\mathbf{s \geqslant 1}$; each $\Theta_{\lambda}$ has a simple top $S_{\lambda}$; and $S_{\lambda} \ncong S_{\mu}$ for $\lambda \neq \mu$.} These conditions are always true for the classical stratifying system of a standardly stratified basic algebra. In particular, in that case $Q = _AA$.

\begin{proposition}
Let $0 \neq M \in \mathcal{F} (\Theta)$ and $i = \text{min (M)}$. Then there is an exact sequence
\begin{equation}
\xymatrix{ 0 \ar[r] & M[1] \ar[r] & M \ar[r] & \bigoplus _{h (\lambda) = i } \Theta_{\lambda} ^{\oplus m_{\lambda}} \ar[r] & 0}
\end{equation}
such that $M[1] \in \mathcal{F} (\Theta)$ and $\text{min} (M[1]) > \text{min} (M)$.
\end{proposition}

\begin{proof}
This is Proposition 2.9 in \cite{Marcos2} which deals with the special case that $\Lambda$ is a linearly ordered set. The general case can be proved similarly by observing the fact that Ext$ _A^1 (\Theta_{\lambda}, \Theta_{\mu}) = 0$ if $h (\lambda) = h (\mu)$.
\end{proof}

It is clear that $m_{\lambda} = [M: \Theta_{\lambda}]$. Based on this proposition, we make the following definition:

\begin{definition}
Let $M \in \mathcal{F} (\Theta)$ with $\text{min} (M) = i$. We say $M$ is generated in height $i$ if in sequence (2.1) we have
\begin{equation*}
\top M = M / \rad M \cong \top \Big{(}\bigoplus _{h (\lambda) = i } \Theta_{\lambda} ^{\oplus m_{\lambda}} \Big{)} = \bigoplus _{h (\lambda) = i } S_{\lambda} ^{\oplus m_{\lambda}}.
\end{equation*}
\end{definition}

We introduce some notation: if $M \in \mathcal{F} (\Theta)$ is generated in height $i$, then define $M_i = \bigoplus _{h (\lambda) = i } \Theta_{\lambda} ^{\oplus m_{\lambda}}$ in sequence (2.1). If $M[1]$ is generated in some height $j$, we can define $M[2] = M[1][1]$ and $M[1]_j$ in a similar way. This procedure can be repeated.

\begin{proposition}
Let $0 \rightarrow L \rightarrow M \rightarrow N \rightarrow 0$ be an exact sequence in $\mathcal{F} (\Theta)$. If $M$ is generated in height $i$, so is $N$. Conversely, if both $L$ and $N$ are generated in height $i$, then $M$ is generated in height $i$ as well.
\end{proposition}

\begin{proof}
We always have $\top N \subseteq \top M$ and $\top M \subseteq \top L \oplus \top N$. The conclusion follows from these inclusions and the rightmost identity in the above definition.
\end{proof}

Notice that $[Q_{\lambda} : \Theta_{\lambda}] = 1$ and $[Q_{\lambda} : \Theta_{\mu}] = 0$ for all $\mu \ngeqslant \lambda$. We claim that $Q_{\lambda}$ is generated in height $h(\lambda)$ for $\lambda \in \Lambda$. Indeed, the algebra $B = \text{End} _A (Q)^{op}$ is a standardly stratified algebra, with projective modules Hom$ _A (Q, Q_{\lambda})$ and standard modules Hom$ _A (Q, \Theta_{\lambda})$, $\lambda \in \Lambda$. Moreover, the functor Hom$ _A (Q, -)$ gives an equivalence between $\mathcal{F} (\Theta) \subseteq A$-mod and $\mathcal{F} (_B \Delta) \subseteq B$-mod. Using this equivalence and the standard filtration structure of projective $B$-modules we deduce the conclusion.

\begin{lemma}
If $M \in \mathcal{F} (\Theta)$ is generated in height $i$ with $[M: \Theta_{\lambda} ] = m_{\lambda}$, then $M$ has a relative projective cover $Q^i \cong \bigoplus _{h (\lambda) = i} Q_{\lambda} ^{\oplus m_{\lambda}}$.
\end{lemma}

\begin{proof}
There is a surjection $f: M \rightarrow \bigoplus _{h(\lambda) = i} \Theta_{\lambda} ^{\oplus m_{\lambda}}$ by Proposition 2.1. Consider the following diagram:
\begin{equation*}
\xymatrix{ & \bigoplus _{h (\lambda) = i} Q_{\lambda} ^{\oplus m_{\lambda}} \ar[d]^p \ar@{-->}[dl]_q \\
M \ar[r]^f & \bigoplus _{h (\lambda) = i} \Theta_{\lambda} ^{\oplus m_{\lambda}} \ar[r] & 0.}
\end{equation*}
Since $Q^i$ is projective in $\mathcal{F} (\Theta)$, the projection $p$ factors through the surjection $f$. In particular, $\top \Big{(} \bigoplus _{h (\lambda) = i} \Theta_{\lambda} ^{\oplus m_{\lambda}} \Big{)} = \bigoplus _{h (\lambda) = i} S_{\lambda} ^{\oplus m_{\lambda}}$ is in the image of $fq$. Since $M$ is generated in height $i$, $f$ induces an isomorphism between $\top M$ and $\top \Big{(} \bigoplus _{h (\lambda) = i} \Theta_{\lambda} ^{\oplus m_{\lambda}} \Big{)}$. Thus $\top M$ is in the image of $q$, and hence $q$ is surjective. It is clear that $q$ is minimal, so $Q^1 = \bigoplus _{h (\lambda) = i} Q_{\lambda} ^{\oplus m_{\lambda}}$ is a relative projective cover of $M$. The uniqueness follows from Proposition 8.3 in \cite{Webb}.
\end{proof}

We use $\Omega _{\Theta} ^i (M)$ to denote the \textit{i-th relative syzygy} of $M$. Actually, for every $M \in \mathcal{F} (\Theta)$ there is always a relative projective cover by Proposition 8.3 in \cite{Webb}.

The following definition is an analogue of \textit{linear modules} in the representation theory of graded algebras.

\begin{definition}
An $A$-module $M \in \mathcal{F} (\Theta)$ is said to be linearly filtered if there is some $i \in \mathbb{N}$ such that $\Omega _{\Theta} ^s (M)$ is generated in height $i + s$ for $s \geqslant 0$.
\end{definition}

Equivalently, $M \in \mathcal{F} (\Theta)$ is linearly filtered if and only if it is generated in height $i$  and has a relative projective resolution
\begin{equation*}
\xymatrix{ 0 \ar[r] & Q^l \ar[r] & Q^{l-1} \ldots \ar[r] & Q^{i+1} \ar[r] & Q^i \ar[r] & M \ar[r] & 0}
\end{equation*}
such that each $Q^s$ is generated in height $s$, $i \leqslant s \leqslant l$.

We remind the reader that there is a common upper bound for the relative projective dimensions of modules contained in $\mathcal{F} (\Theta)$, which is the length of the longest chains in the finite poset $(\Lambda, \leqslant)$. It is also clear that if $M$ is linearly filtered, so are all relative syzygies and direct summands. In other words, the subcategory $\mathcal{LF} (\Theta)$ constituted of linearly filtered modules contains all relative projective modules, and is closed under summands and relative syzygies. But in general it is not closed under extensions, kernels of epimorphisms and cokernels of monomorphisms.

\begin{proposition}
Let $0 \rightarrow L \rightarrow M \rightarrow N \rightarrow 0$ be an exact sequence in $\mathcal{F} (\Theta)$ such that all terms are generated in height $i$. If $L$ is linearly filtered, then $M$ is linearly filtered if and only if $N$ is linearly filtered.
\end{proposition}

\begin{proof}
Let $m_{\lambda} = [M : \Theta_{\lambda}]$, $l_{\lambda} = [L : \Theta_{\lambda}]$ and $n_{\lambda} = [N : \Theta_{\lambda}]$. By the previous lemma, we get the following commutative diagram:

\begin{equation*}
\xymatrix{ & 0 \ar[d] & 0 \ar[d] & 0 \ar[d] & \\
0 \ar[r] & \Omega _{\Theta}(L) \ar[r] \ar[d] & \Omega _{\Theta}(M) \ar[r] \ar[d] & \Omega _{\Theta}(N) \ar[r] \ar[d] & 0 \\
0 \ar[r] & \bigoplus _{h(\lambda) = i} Q_{\lambda}^{\oplus l_{\lambda}} \ar[r] \ar[d] & \bigoplus _{h(\lambda) = i} Q_{\lambda} ^{\oplus m_{\lambda}} \ar[r] \ar[d] & \bigoplus _{h(\lambda) = i} Q_{\lambda} ^{\oplus n_{\lambda}} \ar[r] \ar[d] & 0 \\
0 \ar[r] & L \ar[r] \ar[d] & M \ar[r] \ar[d] & N \ar[r] \ar[d] & 0 \\
 & 0 & 0 & 0 &}
\end{equation*}
Since $\Omega _{\Theta} (L)$ is generated in height $i+1$, by Proposition 2.3, $\Omega _{\Theta} (M)$ is generated in height $i+1$ if and only if $\Omega _{\Theta} (N)$ is generated in height $i+1$. Replacing $L$, $M$ and $N$ by $\Omega _{\Theta} (L)$, $\Omega _{\Theta}(M)$ and $\Omega _{\Theta} (N)$ respectively, we conclude that $\Omega _{\Theta}^2 (M)$ is generated in height $i+2$ if and only if $\Omega _{\Theta}^2 (N)$ is generated in height $i+2$. The conclusion follows from induction.
\end{proof}

\begin{corollary}
Suppose that $M \in \mathcal{F} (\Theta)$ is generated in height $i$ and linearly filtered. If $\bigoplus _{h(\lambda) = i} \Theta_{\lambda}$ is linearly filtered, then $M[1]$ is generated in height $i+1$ and linearly filtered.
\end{corollary}

\begin{proof}
Clearly $\bigoplus _{h(\lambda) = i} \Theta_{\lambda}$ is generated in height $i$. Let $m_{\lambda} = [M : \Theta_{\lambda}]$. Notice that both $M$ and $\bigoplus _{h(\lambda) = i} \Theta_{\lambda} ^{\oplus m_{\lambda}}$ have projective cover $\bigoplus _{h(\lambda) = i} Q_{\lambda} ^{\oplus m_{\lambda}}$. Thus the exact sequence
\begin{equation*}
\xymatrix{ 0 \ar[r] & M[1] \ar[r] & M \ar[r] & \bigoplus _{h(\lambda) = i} \Theta_{\lambda} ^{\oplus m_{\lambda}} \ar[r] & 0}
\end{equation*}
induces the following diagram:
\begin{equation*}
\xymatrix{ & \Omega _{\Theta}(M) \ar@{^{(}->} [r] \ar[d] & \Omega _{\Theta} \Big{(} \bigoplus _{h(\lambda) = i} \Theta_{\lambda} ^{\oplus m_{\lambda}} \Big{)} \ar@{->>}[r] \ar[d] & M[1] \\
& \bigoplus _{h(\lambda) = i} Q_{\lambda} ^{\oplus m_{\lambda}} \ar@{=}[r] \ar[d] & \bigoplus _{h(\lambda) = i} Q_{\lambda} ^{\oplus m_{\lambda}} \ar[d] \\
M[1] \ar@{^{(}->}[r] & M \ar@{->>}[r] & \bigoplus _{h(\lambda) = i} \Theta_{\lambda} ^{\oplus m_{\lambda}}.}
\end{equation*}

Consider the top sequence. Since both $\Omega _{\Theta} (\bigoplus _{h(\lambda) = i} \Theta_{\lambda} ^{\oplus m_{\lambda}})$ and $\Omega _{\Theta} (M)$ are generated in height $i+1$ and linearly filtered, $M[1]$ is also generated in height $i+1$ and linearly filtered by Propositions 2.3 and 2.6.
\end{proof}

These results tell us that linearly filtered modules have properties similar to those of linear modules of graded algebras.

\begin{lemma}
Let $M \in \mathcal{F} (\Theta)$ be generated in height $i$ and $m_{\lambda} = [M : \Theta_{\lambda}]$. If Hom$ _A (Q, \Theta) \cong \text{Hom} _A (\Theta, \Theta)$, then
\begin{equation*}
\text{Hom} _A (M, \Theta) \cong \text{Hom} _A (\bigoplus _{h(\lambda) = i} Q_{\lambda} ^{\oplus m_{\lambda}}, \Theta) \cong \text{Hom} _A (\bigoplus _{h(\lambda) = i} \Theta_{\lambda} ^{\oplus m_{\lambda}}, \Theta).
\end{equation*}
\end{lemma}

\begin{proof}
We claim that Hom$_A (Q, \Theta) \cong \text{Hom} _A (\Theta, \Theta)$ as End$ _A (\Theta)$-modules implies Hom$_A (Q_{\lambda}, \Theta) \cong \text{Hom} _A (\Theta_{\lambda}, \Theta)$ for every $\lambda \in \Lambda$. Then the second isomorphism follows immediately. First, notice that Hom$ _A (\Theta, \Theta)$ is a basic algebra with $n$ non-zero indecomposable summands Hom$ _A (\Theta_{\lambda}, \Theta)$, $\lambda \in \Lambda$, where $n$ is the cardinal number of $\Lambda$. But Hom$ _A (Q, \Theta) \cong \bigoplus _{\lambda \in \Lambda} \text{Hom} _A (Q_{\lambda}, \Theta)$ has at least $n$ non-zero indecomposable summands. If Hom$_A (Q, \Theta) \cong \text{Hom} _A (\Theta, \Theta)$, by Krull-Schmidt theorem, Hom$ _A (Q_{\lambda}, \Theta)$ must be indecomposable, and is isomorphic to some Hom$ _A (\Theta_{\mu}, \Theta)$.

If $\lambda \in \Lambda$ is maximal, $\text{Hom} _A (\Theta_{\lambda}, \Theta) \cong \text{Hom} _A (Q_{\lambda}, \Theta)$ since $Q_{\lambda} \cong \Theta_{\lambda}$. Define $\Lambda_1$ to be the subset of maximal elements in $\Lambda$ and consider $\lambda_1 \in \Lambda \setminus \Lambda_1$ which is maximal. We have
\begin{equation*}
\text{Hom} _A (Q_{\lambda_1}, \Theta) \ncong \text{Hom} _A (\Theta_{\lambda}, \Theta) \cong \text{Hom} _A (Q_{\lambda}, \Theta)
\end{equation*}
for every $\lambda \in \Lambda_1$ since End$ _A (\Theta)$ is a basic algebra. Therefore, $\text{Hom} _A (Q_{\lambda_1}, \Theta)$ must be isomorphic to some $\text{Hom} _A (\Theta_{\mu}, \Theta)$ with $\mu \in \Lambda \setminus \Lambda_1$. But $\text{Hom} _A (\Theta_{\mu}, \Theta)$ contains a surjection from $\Theta_{\mu}$ to the direct summand $\Theta_{\mu}$ of $\Theta$, and Hom$ _A (Q _{\lambda_1}, \Theta)$ contains a surjection from $Q_{\lambda_1}$ to $\Theta_{\mu}$ if and only if $\lambda_1 = \mu$. Thus we get $\lambda_1 = \mu$. Repeating the above process, we have Hom$_A (Q_{\lambda}, \Theta) \cong \text{Hom} _A (\Theta_{\lambda}, \Theta)$ for every $\lambda \in \Lambda$.

Applying Hom$ _A (-, \Theta)$ to the surjection $M \rightarrow \bigoplus _{h(\lambda) = i} \Theta_{\lambda} ^{\oplus m_{\lambda}}$ we get
\begin{equation*}
\text{Hom} _A (\bigoplus _{h(\lambda) = i} \Theta_{\lambda} ^{\oplus m_{\lambda}}, \Theta) \subseteq \text{Hom} _A (M, \Theta).
\end{equation*}
Similarly, from the relative projective covering map $\bigoplus _{h(\lambda) = i} Q_{\lambda} ^{\oplus m_{\lambda}} \rightarrow M$ we have
\begin{equation*}
\text{Hom} _A (M, \Theta) \subseteq \text{Hom} _A (\bigoplus _{h(\lambda) = i} Q_{\lambda} ^{\oplus m_{\lambda}}, \Theta).
\end{equation*}
Comparing these two inclusions and using the second isomorphism, we deduce the first isomorphism.
\end{proof}

The reader may aware that the above lemma is an analogue to the following result in representation theory of graded algebras: if $A$ is a graded algebra and $M$ is a graded module generated in degree 0, then Hom$ _A (M, A_0) \cong \text{Hom} _A (M_0, A_0)$.

\begin{lemma}
Suppose that Hom$ _A (Q, \Theta) \cong \text{Hom} _A (\Theta, \Theta)$. If $M \in \mathcal{F} (\Theta)$ is generated in height $i$, then Ext$ _A^s (M, \Theta) \cong \text{Ext} _A^{s-1} (\Omega _{\Theta}(M), \Theta)$ for all $s \geqslant 1$.
\end{lemma}

\begin{proof}
Let $m_{\lambda} = [M: \Theta_{\lambda}]$. Applying Hom$ _A (-, \Theta)$ to the exact sequence
\begin{equation*}
\xymatrix{ 0 \ar[r] & \Omega _{\Theta} (M) \ar[r] & \bigoplus _{h(\lambda) = i} Q_{\lambda} ^{\oplus m_{\lambda}} \ar[r] & M \ar[r] & 0}
\end{equation*}
we get a long exact sequence. In particular, for all $s \geqslant 2$, by observing the segment
\begin{align*}
& 0 = \text{Ext} _A^{s-1} (\bigoplus _{h(\lambda) = i} Q ^{\oplus m_{\lambda}}, \Theta) \rightarrow \text{Ext} _A^{s-1} (\Omega _{\Theta} (M), \Theta) \\
& \rightarrow \text{Ext} _A^s (M, \Theta) \rightarrow \text{Ext} _A^s (\bigoplus _{h(\lambda) = i} Q_{\lambda} ^{\oplus m_{\lambda}}, \Theta) = 0
\end{align*}
we conclude $\text{Ext} _A^{s-1} (\Omega _{\Theta}(M), \Theta) \cong \text{Ext} _A^s (M, \Theta)$.

For $s=1$, we have
\begin{align*}
& 0 \rightarrow \text{Hom} _A (M, \Theta) \rightarrow \text{Hom} _A (\bigoplus _{h(\lambda) = i} Q_{\lambda} ^{\oplus m_{\lambda}}, \Theta)\\
& \rightarrow \text{Hom} _A (\Omega _{\Theta} (M), \Theta) \rightarrow \text{Ext} _A^1 (M, \Theta) \rightarrow 0.
\end{align*}
By the previous lemma, the first inclusion is an ismorphism. Thus Ext$ _A^1 (M, \Theta) \cong \text{Hom} _A (\Omega _{\Theta} (M), \Theta)$.
\end{proof}

We remind the reader that although $\Gamma = \text{Ext} _A ^{\ast} (\Theta, \Theta)$ has a natural grading, the classical Koszul theory cannot be applied directly since $\Gamma_0 = \text{End} _A (\Theta)$ may not be a semisimple algebra. Thus we introduce \textit{generalized Koszul algebras} as follows:

\begin{definition}
Let $R = \bigoplus _{i \geqslant 0} R_i$ be a positively graded locally finite $k$-algebra, i.e., $\dim _k R_i < \infty$ for each $i \geqslant 0$. A graded $R$-module $M$ is said to be linear if it has a linear projective resolution
\begin{equation*}
\xymatrix{ \ldots \ar[r] & P^s \ar[r] & \ldots \ar[r] & P^1 \ar[r] & P^0 \ar[r] & M}
\end{equation*}
such that $P^s$ is generated in degree $s$. The algebra $R$ is said to be generalized Koszul if $R_0$ viewed as a $R$-module has a linear projective resolution.
\end{definition}

It is easy to see from this definition that $M$ is linear if and only if $\Omega^s (M)$ is generated in degree $s$ for all $s \geqslant 0$.

\begin{proposition}
Suppose that Hom$ _A (Q, \Theta) \cong \text{Hom} _A (\Theta, \Theta)$, and $\Theta_{\lambda}$ are linearly filtered for all $\lambda \in \Lambda$. If $M \in \mathcal{F} (\Theta)$ is linearly filtered, then
\begin{equation*}
\text{Ext} _A^{i+1} (M, \Theta) = \text{Ext} _A^1 (\Theta, \Theta) \cdot \text{Ext} _A^i (M, \Theta)
\end{equation*}
for all $i \geqslant 0$, i.e., $\text{Ext} _A^{\ast} (M, \Theta)$ as a graded $\Gamma = \text{Ext} _A^{\ast} (\Theta, \Theta)$-module is generated in degree 0.
\end{proposition}

\begin{proof}
Suppose that $M$ is generated in height $d$ and linearly filtered. By Lemma 2.9,
\begin{equation*}
\text{Ext} _A^{i+1} (M, \Theta) \cong \text{Ext} _A^i (\Omega _{\Theta} (M), \Theta).
\end{equation*}
But $\Omega _{\Theta}$ is generated in height $d+1$ and linearly filtered. Thus by induction
\begin{equation*}
\text{Ext} _A^{i+1} (M, \Theta) \cong \text{Ext} _A^1 (\Omega^i _{\Theta} (M), \Theta), \quad \text{Ext} _A^i (M, \Theta) \cong \text{Hom} _A (\Omega^i _{\Theta} (M), \Theta).
\end{equation*}
Therefore, it suffices to show
\begin{equation*}
\text{Ext} _A^1 (M, \Theta) = \text{Ext} _A^1 (\Theta, \Theta) \cdot \text{Hom} _A (M, \Theta)
\end{equation*}
since we can replace $M$ by $\Omega _{\Theta} ^i (M)$, which is linearly filtered as well.\

Let $m_{\lambda} = [M: \Theta_{\lambda}]$ and define $Q^0 = \bigoplus _{h(\lambda) = d} Q_{\lambda} ^{\oplus m_{\lambda}}$, $M_0 = \bigoplus _{h(\lambda) = d} \Theta_{\lambda} ^{\oplus m_{\lambda}}$. We have the following commutative diagram:
\begin{equation*}
\xymatrix {0 \ar[r] & \Omega _{\Theta} (M) \ar@{=}[d] \ar[r] & Q^0[1] \ar[r] \ar[d] & M[1] \ar[d] \ar[r] & 0\\
0 \ar[r] & \Omega _{\Theta} (M) \ar[r] & Q^0 \ar[r] \ar[d] & M \ar[d] \ar[r] & 0\\
 & & M_0 \ar@{=}[r] & M_0}
\end{equation*}
where $Q^0[1] = \Omega _{\Theta} (M_0)$, see Proposition 2.1.

Observe that all terms in the top sequence are generated in height $d+1$ and linearly filtered. For every $\lambda \in \Lambda$ with $h(\lambda) = d+1$, we have
\begin{equation*}
[\Omega _{\Theta} (M): \Theta_{\lambda}] + [M[1]: \Theta_{\lambda}] = [Q^0[1]: \Theta_{\lambda}].
\end{equation*}
Let $r_{\lambda}, s_{\lambda}$ and $t_{\lambda}$ be the corresponding numbers in the last equality. Then we get a split short exact sequence
\begin{equation*}
\xymatrix{ 0 \ar[r] & \bigoplus _{h(\lambda) = d+1} \Theta_{\lambda} ^{\oplus r_{\lambda}} \ar[r] & \bigoplus _{h(\lambda) = d+1} \Theta_{\lambda} ^{\oplus t_{\lambda}} \ar[r] & \bigoplus _{h(\lambda) = d+1} \Theta_{\lambda} ^{\oplus s_{\lambda}} \ar[r] & 0}.
\end{equation*}
Applying Hom$ _A (-, \Theta)$ to this sequence and using Lemma 2.8, we obtain the exact sequence
\begin{equation*}
0 \rightarrow \text{Hom} _A (M[1], \Theta) \rightarrow \text{Hom} _A (Q^0[1], \Theta) \rightarrow \text{Hom} _A (\Omega _{\Theta} (M), \Theta) \rightarrow 0.
\end{equation*}
Therefore, each map $\Omega _{\Theta} (M) \rightarrow \Theta$ can extend to a map $Q^0[1] \rightarrow \Theta$.

To prove Ext$_A^1 (M, \Theta) = \text{Ext} _A^1 (\Theta, \Theta) \cdot \text{Hom} _A (M, \Theta)$, by Lemma 2.9 we first identify $\text{Ext} _A^1 (M, \Theta)$ with $\text{Hom} _A (\Omega _{\Theta} (M), \Theta)$. Take an element $x \in \text{Ext} _A^1 (M, \Theta)$ and let $g: \Omega _{\Theta} (M) \rightarrow \Theta$ be the corresponding homomorphism. As we just showed, $g$ can extend to $Q^0[1]$, and hence there is a homomorphism $\tilde{g}: Q^0[1] \rightarrow \Theta$ such that $g = \tilde{g} \iota$, where $\iota$ is the inclusion.
\begin{align*}
\xymatrix { \Omega _{\Theta} (M) \ar[r]^{\iota} \ar[d]^g & Q^0[1] \ar[dl]^{\tilde{g}}\\
\Theta.}
\end{align*}

We have the following commutative diagram:
\begin{align*}
\xymatrix{ 0 \ar[r] & \Omega _{\Theta} (M) \ar[d]^{\iota} \ar[r] & Q^0 \ar@{=}[d] \ar[r] & M \ar[d]^p \ar[r] & 0\\
0 \ar[r] & Q^0[1] \ar[r] & Q^0 \ar[r] & M_0 \ar[r] & 0,}
\end{align*}
where $p$ is the projection of $M$ onto $M_0$. The map $\tilde{g}: Q^0[1] \rightarrow \Theta$ gives a push-out of the bottom sequence:
\begin{align*}
\xymatrix{ 0 \ar[r] & \Omega _{\Theta} (M) \ar[d]^{\iota} \ar[r] & Q^0 \ar@{=}[d] \ar[r] & M \ar[d]^p \ar[r] & 0\\
0 \ar[r] & Q^0[1] \ar[r] \ar[d] ^{\tilde{g}} & Q^0 \ar[d] \ar[r] & M_0 \ar[r] \ar@{=}[d] & 0\\
0 \ar[r] & \Theta \ar[r] & E \ar[r] & M_0 \ar[r] & 0.}
\end{align*}

Since $M_0 \cong \bigoplus _{h(\lambda) = d} \Theta_{\lambda} ^{\oplus m_{\lambda}}$, the bottom sequence represents some
\begin{equation*}
y \in \text{Ext} _A^1 (\Theta ^{\oplus m}, \Theta) = \bigoplus _{i=1} ^m \text{Ext} _A^1 (\Theta, \Theta)
\end{equation*}
where $m = \sum _{h(\lambda) = d} m_{\lambda}$. Therefore, we can write $y = y_1 + \ldots + y_m$ where $y_i \in \text{Ext} _A^1 (\Theta, \Theta)$ is represented by the sequence
\begin{equation*}
\xymatrix{ 0 \ar[r] & \Theta \ar[r] & E_i \ar[r] & \Theta \ar[r] & 0}.
\end{equation*}
Composed with the inclusions $\epsilon_{\lambda}: \Theta_{\lambda} \rightarrow \Theta$, we get the map $(p_1, \ldots, p_m)$ where each component $p_i$ is defined in an obvious way. Consider the pull-backs:
\begin{equation*}
\xymatrix {0 \ar[r] & \Theta \ar[r] \ar@{=}[d] & F_i \ar[r] \ar[d] & M \ar[r] \ar[d]^{p_i} & 0 \\
0 \ar[r] & \Theta \ar[r] & E_i \ar[r] & \Theta \ar[r] & 0.}
\end{equation*}
Denote by $x_i$ the top sequence. Then
\begin{equation*}
x = \sum_{i=1}^m x_i = \sum_{i=1}^m y_i \cdot p_i \in \text{Ext} _A^1 (\Theta, \Theta) \cdot \text{Hom} _A (M, \Theta),
\end{equation*}
so $\text{Ext} _A^1 (M, \Theta) \subseteq \text{Ext} _A^1 (\Theta, \Theta) \cdot \text{Hom} _A (M, \Theta)$. The other inclusion is obvious.
\end{proof}

Now we can prove the main result.

\begin{theorem}
Let $(\underline {\Theta}, \underline{Q})$ be an EPSS indexed by a finite poset $(\Lambda, \leqslant)$ such that Ext$_A ^i (Q, \Theta) =0$ for all $i \geqslant 1$ and Hom$ _A (Q, \Theta) \cong \text{Hom} _A (\Theta, \Theta)$. Suppose that all $\Theta_{\lambda}$ are linearly filtered for $\lambda \in \Lambda$. If $M \in \mathcal{F} (\Theta)$ is linearly filtered, then the graded module $\text{Ext} _A^{\ast} (M, \Theta)$ has a linear projective resolution. In particular, $\Gamma = \text{Ext} _A^{\ast} (\Theta, \Theta)$ is a generalized Koszul algebra.
\end{theorem}

\begin{proof}
Suppose that $M$ is generated in height $d$. Define $m_{\lambda} = [M: \Theta_{\lambda}]$ for $\lambda \in \Lambda$, $Q^0 = \bigoplus _{h(\lambda) = d} Q_{\lambda} ^{\oplus m_{\lambda}}$, and $M_0 = \bigoplus _{h(\lambda) = d} \Theta_{\lambda} ^{\oplus m_{\lambda}}$. As in the proof of the previous lemma, we have the following short exact sequence of linearly filtered modules generated in height $d+1$:
\begin{align*}
\xymatrix{ 0 \ar[r] & \Omega _{\Theta} (M) \ar[r] & \Omega _{\Theta} (M_0) \ar[r] & M[1] \ar[r] & 0}
\end{align*}
where $\Omega _{\Theta} (M_0) = Q^0[1]$. This sequence induces exact sequences recursively (see the proof of Proposition 2.6):
\begin{align*}
\xymatrix{ 0 \ar[r] & \Omega^i _{\Theta} (M) \ar[r] & \Omega^i _{\Theta} (M_0) \ar[r] & \Omega^{i-1} _{\Theta} (M[1]) \ar[r] & 0,}
\end{align*}
where all modules are linearly filtered and generated in height $d+i$. Again as in the proof of the previous lemma, we get an exact sequence
\begin{align*}
0 \rightarrow \text{Hom} _A (\Omega _{\Theta} ^{i-1} (M[1]), \Theta) \rightarrow \text{Hom} _A (\Omega _{\Theta} ^i (M_0), \Theta) \rightarrow \text{Hom} _A (\Omega _{\Theta} ^i (M), \Theta) \rightarrow 0.
\end{align*}
According to Lemma 2.9, the above sequence is isomorphic to:
\begin{align*}
0 \rightarrow \text{Ext} ^{i-1} _A (M[1], \Theta) \rightarrow \text{Ext}^i _A (M_0, \Theta) \rightarrow \text{Ext}^i _A (M, \Theta) \rightarrow 0.
\end{align*}
Now let the index $i$ vary and put these sequences together. We have:
\begin{align*}
\xymatrix{ 0 \ar[r] & E(M[1]) \langle 1 \rangle \ar[r] & E(M_0) \ar[r]^p & E(M) \ar[r] & 0,}
\end{align*}
where $E = \text{Ext} _A^{\ast} (-, \Theta)$ and $\langle - \rangle$ is the degree shift functor of graded modules. That is, for a graded module $T = \bigoplus _{i \geqslant 0} T_i$, $T \langle 1 \rangle _i$ is defined to be $T_{i-1}$.

Since $M_0 \in \text{add} (\Theta)$, $E(M_0)$ is a projective $\Gamma$-module. It is generated in degree 0 by the previous lemma. Similarly, $E(M[1])$ is generated in degree 0, so $E(M[1]) \langle 1 \rangle$ is generated in degree 1. Therefore, the map $p$ is a graded projective covering map. Consequently, $\Omega(E(M)) \cong E(M[1]) \langle 1 \rangle$ is generated in degree 1.

Replacing $M$ by $M[1]$ (since it is also linearly filtered), we have
\begin{equation*}
\Omega^2 (E(M)) \cong \Omega (E(M[1]) \langle 1 \rangle) \cong \Omega( E(M[1]) \langle 1 \rangle \cong E(M[2]) \langle 2 \rangle,
\end{equation*}
which is generated in degree 2. By recursion, $\Omega^i (E(M)) \cong E(M [i]) \langle i \rangle$ is generated in degree $i$ for all $i \geqslant 0$. Thus $E(M)$ is a linear $\Gamma$-module.

In particular let $M = Q_{\lambda}$ for a certain $\lambda \in \Lambda$. We get
\begin{equation*}
E(Q_{\lambda}) = \text{Ext} _A^{\ast} (Q_{\lambda}, \Theta) = \text{Hom} _A (Q_{\lambda}, \Theta)
\end{equation*}
is a linear $\Gamma$-module. Therefore,
\begin{align*}
\bigoplus _{\lambda \in \Lambda} E (Q_{\lambda}, \Theta) & = \bigoplus _{\lambda \in \Lambda} \text{Hom} _A (Q_{\lambda}, \Theta) \cong \text{Hom} _A (\bigoplus _{\lambda \in \Lambda} Q_{\lambda}, \Theta)\\
& = \text{Hom} _A (Q, \Theta) \cong \text{Hom} _A (\Theta, \Theta) = \Gamma_0
\end{align*}
is a linear $\Gamma$-module. So $\Gamma$ is a generalized Koszul algebra.
\end{proof}

\begin{remark}
To get the above result we made some assumptions on the EPSS $(\underline {\Theta}, \underline{Q})$. Firstly, each $\Theta_{\lambda}$ has a simple top $S_{\lambda}$ and $S_{\lambda} \ncong S_{\mu}$ for $\lambda \neq \mu$; secondly, Ext$ _A^s (Q, \Theta) = 0$ for every $s \geqslant 1$. These two conditions always hold for standardly stratified basic algebras. We also suppose that Hom$ _A (\Theta, \Theta) \cong \text{Hom} _A (Q, \Theta)$. This may not be true even if $A$ is a quasi-hereditary algebra.
\end{remark}

Although $\Gamma$ is proved to be a generalized Koszul algebra, in general it does not have the Koszul duality. Consider the following example:

\begin{example}
Let $A$ be the path algebra of the following quiver with relation $\alpha \cdot \beta = 0$. Put an order $x < y < z$.
\begin{equation*}
\xymatrix{ x \ar@/^/ [rr]^{\alpha} & & y \ar@ /^/ [ll]^{\beta} \ar[r] ^{\gamma} & z.}
\end{equation*}
The projective modules and standard modules of $A$ are described as follows:
\begin{equation*}
P_x = \begin{matrix}  & x & \\ & y & \\ x & & z \end{matrix} \qquad P_y = \begin{matrix} & y & \\ x & & z \end{matrix} \qquad P_z = z
\end{equation*}
\begin{equation*}
\Delta_x = x \qquad \Delta_y = \begin{matrix} y \\ x \end{matrix} \qquad \Delta_z = z \cong P_z.
\end{equation*}
This algebra is quasi-hereditary. Moreover, Hom$ _A (\Delta, \Delta) \cong \text{Hom} _A (A, \Delta)$, and all standard modules are linearly filtered. Therefore, $\Gamma = \text{Ext} _A^{\ast} (\Delta, \Delta)$ is a generalized Koszul algebra by the previous theorem.\

We explicitly compute the extension algebra $\Gamma$. It is the path algebra of the following quiver with relation $\gamma \cdot \alpha = 0$.
\begin{equation*}
\xymatrix{ x \ar@/^/ [rr]^{\alpha} \ar@ /_/ [rr]_{\beta} & & y \ar[r] ^{\gamma} & z.}
\end{equation*}
\begin{equation*}
_{\Gamma} P_x = \begin{matrix}  & x_0 & \\ y_0 & & y_1 \\ z_1 & & \end{matrix} \qquad _{\Gamma} P_y = \begin{matrix} y_0 \\ z_1 \end{matrix} \qquad _{\Gamma} P_z = z_0
\end{equation*}
and
\begin{equation*}
_{\Gamma} \Delta_x = x_0 \qquad _{\Gamma} \Delta_y = y_0 \qquad _{\Gamma} \Delta_z = z_0 \qquad \Gamma_0 = \begin{matrix} x_0 \\ y_0 \end{matrix} \oplus y_0 \oplus z_0 \ncong _{\Gamma} \Delta.
\end{equation*}
Here we use indices to mark the degrees of simple composition factors. As asserted by the theorem, $\Gamma_0$ has a linear projective resolution. But $_{\Gamma} \Delta$ is not a linear $\Gamma$-module (we remind the reader that the two simple modules $y$ appearing in $_{\Gamma} P_x$ lie in different degrees!).

By computation, we get the extension algebra $\Gamma' = \text{Ext} _{\Gamma} ^{\ast} (\Gamma_0, \Gamma_0)$, which is the path algebra of the following quiver with relation $\beta \cdot \alpha =0$.
\begin{equation*}
\xymatrix {x \ar[r] ^{\alpha} & y \ar[r] ^{\beta} & z.}
\end{equation*}
Since $\Gamma'$ is a Koszul algebra in the classical sense, the Koszul duality holds in $\Gamma'$. It is obvious that the Koszul dual algebra of $\Gamma'$ is not isomorphic to $\Gamma$. Therefore, as we claimed, the Koszul duality does not hold in $\Gamma$.
\end{example}

Let us return to the question of whether $\Gamma = \text{Ext} _A ^{\ast} (\Theta, \Theta)$ is standardly stratified with respect to $\leqslant$. According to Proposition 1.3, this happens if and only if for each pair $\Theta_{\lambda}, \Theta_{\mu}$ and $s \geqslant 0$, Ext$ _A^s (\Theta_{\lambda}, \Theta_{\mu})$ is a projective End$ _A (\Theta_{\mu})$-module. Putting direct summands together, we conclude that $\Gamma$ is standardly stratified with respect to $\leqslant$ if and only if Ext$ _A^s (\Theta, \Theta)$ is a projective $\bigoplus _{\lambda \in \Lambda} \text{End} _A (\Theta_{\lambda})$-module. With the conditions in Theorem 2.12, Ext$ _A^s (\Theta_{\lambda}, \Theta) \cong \text{Hom} _A (\Omega ^s _{\Theta _{\lambda}} (\Theta), \Theta)$ for all $s \geqslant 0$ and $\lambda \in \Lambda$ by Lemma 2.9. Notice that $\Omega ^s _{\Theta} (\Theta _{\lambda})$ is linearly filtered. Suppose that min$(\Omega ^s _{\Theta} (\Theta _{\lambda}) ) = d$ and $m_{\mu} = [\Omega ^s _{\Theta} (\Theta _{\lambda}): \Theta_{\mu}]$. Then
\begin{equation}
\text{Ext} _A^s (\Theta _{\lambda}, \Theta) \cong \text{Hom} _A (\Omega ^s _{\Theta} (\Theta _{\lambda}), \Theta) \cong \text{Hom} _A ( \bigoplus _{h(\mu) = d} \Theta_{\mu} ^{\oplus m_{\mu}}, \Theta),
\end{equation}
which is a projective $\Gamma_0 = \text{End}_A (\Theta)$-module.

With this observation, we have:

\begin{corollary}
Let $(\underline {\Theta}, \underline {Q})$ be an EPSS indexed by a finite poset $(\Lambda, \leqslant)$. Suppose that all $\Theta_{\lambda}$ are linearly filtered for $\lambda \in \Lambda$, and Hom$ _A (Q, \Theta) \cong \text{Hom} _A (\Theta, \Theta)$. Then $\Gamma = \text{Ext} _A^{\ast} (\Theta, \Theta)$ is standardly stratified for $\leqslant$ if and only if End$ _A (\Theta)$ is a projective $\bigoplus _{\lambda \in \Lambda} \text{End} _A (\Theta_{\lambda})$-module.
\end{corollary}

\begin{proof}
If $\Gamma$ is standardly stratified for $\leqslant$, then in particular $\Gamma_0 = \text{End} _A (\Theta)$ is a projective $\bigoplus _{\lambda \in \Lambda} \text{End} _A (\Theta_{\lambda})$-module by Proposition 1.3. Conversely, if $\Gamma_0 = \text{End} _A (\Theta)$ is a projective $\bigoplus _{\lambda \in \Lambda} \text{End} _A (\Theta_{\lambda})$-module, then by the isomorphism in (2.2) Ext$ _A^s (\Theta, \Theta) = \bigoplus _{\lambda \in \Lambda} \text{Ext} _A^s (\Theta_{\lambda}, \Theta)$ is a projective $\Gamma_0$-module for all $s \geqslant 0$, so it is a projective $\bigoplus _{\lambda \in \Lambda} \text{End} _A (\Theta_{\lambda})$-module as well. Again by Proposition 1.3, $\Gamma$ is standardly stratified with respect to $\leqslant$.
\end{proof}

If $A$ is quasi-hereditary with respect to $\leqslant$ such that all standard module are linearly filtered, then $\Gamma = \text{Ext} _A^{\ast} (\Delta, \Delta)$ is again quasi-hereditary for this partial order by Corollary 1.5, and $\Gamma_0$ has a linear projective resolution by the previous theorem. Let $_{\Gamma} \Delta$ be the direct sum of all standard modules of $\Gamma$ with respect to $\leqslant$. The reader may wonder whether $_{\Gamma} \Delta$ has a linear projective resolution as well. The following proposition gives a partial answer to this question.

\begin{proposition}
With the above notation, if $_{\Gamma} \Delta$ has a linear projective resolution, then $\Gamma_0 \cong  _{\Gamma}\Delta$, or equivalently Hom$ _A (\Delta_{\lambda}, \Delta_{\mu}) \neq 0$ only if $\lambda = \mu$, $\lambda, \mu \in \Lambda$. If furthermore Hom$ _A (A, \Delta) \cong \text{End} _A (\Delta)$, then $\Delta \cong A / \rad A$.
\end{proposition}

\begin{proof}
We have proved that the $k$-linear category associated to $\Gamma$ is directed with respect to $\leqslant$. By Proposition 1.3, standard modules of $\Gamma$ for $\leqslant$ are exactly indecomposable summands of $\bigoplus _ {\lambda \in \Lambda} \text{End} _A (\Delta_{\lambda})$, i.e., $_{\Gamma} \Delta \cong \bigoplus _ {\lambda \in \Lambda} \text{End} _A (\Delta_{\lambda}) \cong \bigoplus _{\lambda \in \Lambda} k_{\lambda}$. Clearly, $_{\Gamma} \Delta \subseteq \Gamma_0 = \text{End} _A (\Delta)$. If $_{\Gamma} \Delta$ has a linear projective resolution, then by Corollary 2.4 and Remark 2.7 in \cite{Li1}, $_{\Gamma} \Delta$ is a projective $\Gamma_0$-module. Consequently, every summand $k_{\lambda}$ is a projective $\Gamma_0$-module. Since both $_{\Gamma} \Delta$ and $\Gamma_0$ have exactly $| \Lambda |$ pairwise non-isomorphic indecomposable summands, we deduce $_{\Gamma} \Delta \cong \Gamma_0 \cong \bigoplus _{\lambda \in \Lambda} k_{\lambda}$, or equivalently $\text{Hom} _A (\Delta_{\lambda}, \Delta_{\mu}) = 0$ if $\lambda \neq \mu$.

If furthermore Hom$ _A (A, \Delta) \cong \text{End} _A (\Delta)$, then
\begin{equation*}
\Delta \cong \text{Hom} _A (A, \Delta) \cong \text{End} _A (\Delta) \cong \bigoplus _{\lambda \in \Lambda} k_{\lambda} \cong A / \rad A.
\end{equation*}
\end{proof}


\begin{thebibliography}{99}
\bibitem{Abe} N. Abe, \textit{First extension groups of Verma modules and $R$-polynomials}, ArXiv:1003.0169;\
\bibitem{Agoston1} I. \'{A}goston, V. Dlab, and E. Luk\'{a}s, \textit{Quasi-hereditary Extension Algebras}, Alg. Repn. Theorey 6 (2003), 97-117;\
\bibitem{Agoston2} I. \'{A}goston, V. Dlab, and E. Luk\'{a}s, \textit{Standardly stratified extension algebras}, Comm. Algebra 33 (2005), 1357-1368;\
\bibitem{Cline} E. Cline, B. Parshall, and L. Scott, \textit{Stratifying Endomorphism Algebras}, Mem. Amer. Math. Sco. 124 (1996), no. 591;\
\bibitem{Dlab} V. Dlab, and C. Ringel, \textit{The module theoretical approach to Quasi-hereditary Algebras}, Representations of algebras and related topics (Kyoto, 1990), 200-224, London Math. Soc. Lecture Note Ser. 168, Cambridge Univ. Press 1992;\
\bibitem{Drozd} Y. Drozd, V. Mazorchuk, \textit{Koszul duality for extension algebras of standard modules}, J. Pure Appl. Algebra 211 (2007), 484-496;\
\bibitem{Erdmann} K. Erdmann, and C. S\'{a}enz, \textit{On standardly stratified algebras}, Comm. Algebras 31 (2003), 3429-3446;\
\bibitem{Frisk1} A. Frisk, \textit{Two-step tilting for standardly stratified algebras}, Algebra Discrete Math. 2004, no. 3, 38-59;\
\bibitem{Frisk2} A. Frisk, \textit{Dlab's theorem and tilting modules for stratified algebras}, J. Algebra 314 (2007), 507-537;\
\bibitem{Green} E. L. Green, I. Reiten, and {\O}. Solberg, \textit{Dualities on Generalized Koszul Algebras}, Mem. Amer. Math. Soc. 159 (2002), xvi+67pp;\
\bibitem{Klamt} A. Klamt, C. Stroppel, \textit{On the Ext algebras of parabolic Verma modules and $A_{\infty}$-structures}, J. Pure Appl. Algebra 216 (2012), 323-336;\
\bibitem{Li1} L. Li, \textit{A generalized Koszul theory and its application}, ArXiv: 1109.5760;\
\bibitem{Li2} L. Li, \textit{Algebras stratified for all preorders}, accepted by Alg. Repn. Theory, ArXiv: 1110.6501;\
\bibitem{Madsen1} D. Madsen, \textit{Ext-algebras and derived equivalences}, Colloq. Math. 104 (2006), 113-140;\
\bibitem{Madsen2} D. Madsen, \textit{On a Common Generalization of Koszul Duality and Tilting Equivalence}, Adv. Math. 227 (2011), 2327-2348;\
\bibitem{Madsen3} D. Madsen, \textit{Quasi-hereditary algebras and generalized Koszul duality}, ArXiv: 1201.0441;\
\bibitem{Marcos1} E. Marcos, O. Mendoza, and C. S\'{a}enz, \textit{Stratifying systems via relative simple modules}, J. Algebra 280 (2004), 472-487;\
\bibitem{Marcos2} E. Marcos, O. Mendoza, and C. S\'{a}enz, \textit{Stratifying systems via relative projective modules}, Comm. Algebra 33 (2005), 1559-1573;\
\bibitem{Mazorchuk1} V. Mazorchuk, \textit{Some homological properties of the category $\mathcal{O}$}, Pacific J. Math. 232 (2007), 313-341;\
\bibitem{Mazorchuk2} V. Mazorchuk, S. Ovsienko, and C. Stroppel, \textit{Quadratic duals, Koszul dual functors, and applications}, Trans. Amer. Math. Soc. 361 (2009), 1129-1172;\
\bibitem{Mazorchuk3} V. Mazorchuk, \textit{Koszul duality for stratified algebras I: balanced quasi-hereditary algebras}, Manuscripta Math. 131 (2010), 1-10;\
\bibitem{Mazorchuk4} V.Mazorchuk, \textit{Koszul duality for stratified algebras II: standardly stratified algebras}, J. Aust. Math. Soc. 89 (2010),
23-49;\
\bibitem{Miemietz} V. Miemietz, W. Turner, \textit{The Weyl extension algebra of $GL_2(\bar{\mathbb{F}}_p)$}, ArXiv: 1106.5665;
\bibitem{Webb} P. Webb, \textit{Stratifications and Mackey functors. I: functors for a single group}, Proc. London Math. Soc. (3) 82 (2001), 299-336.
\end{thebibliography}
\end{document}